\documentclass[12pt]{article}

\usepackage{amsmath,amssymb,latexsym,theorem,epsfig,amscd}
\usepackage[all]{xy}
\CompileMatrices

\newcommand{\be}{\begin{equation}}
\newcommand{\ee}{\end{equation}}
\newcommand{\eq}[1]{Eq. (\ref{eq-#1})}
\newcommand{\bbc}{{\mathbb C}}
\newcommand{\bbz}{{\mathbb Z}}

\newcommand{\bbr}{{\mathbb R}}
\newcommand{\rk}{{\text{ rk} }}
\newcommand{\Tr}{{\text{ Tr}~ }}
\renewcommand{\Re}{{\text{ Re}~ }}
\newcommand{\End}{{\text{ End}~ }}
\newcommand{\Coh}{{\text{ Coh} }}

\newcommand{\ch}{{\text{ch} }}

\newcommand{\CH}{{\it CH}}
\newcommand{\CHu}{{\it \overline{CH}}}
\newcommand{\CHi}{{\it\underline{CH}}}
\newcommand{\PCH}{{\it PCH}}
\newcommand{\PCHu}{{\it \overline{PCH}}}
\newcommand{\PCHi}{{\it\underline{PCH}}}
\newcommand{\ATT}{{\it ATT}}
\newcommand{\KC}{{\it KC}}
\newcommand{\tH}{{\tilde H}}

\newcommand{\nch}{{\text{ch}}}
\newcommand{\td}{{\text{td}\; }}

\newcommand{\cp}[1]{{\mathbb P}^{#1}}

\newcommand{\om}{{\mathcal O}_M}

\newcommand{\cH}{{\mathcal H}}

\newcommand\Hom{{\rm Hom}}
\newcommand\Ext{{\rm Ext}}

\newcommand{\les}[8]{\xymatrix{       &      & ...  \ar[r]  &  {#1}    \ar@{->} `r[d] `[l] `^dl[dlll]  `^dr/14pt[dll]    [dll] \\
&  {#2} \ar[r] & {#3} \ar[r] & {#4}  \ar `r/10pt[d] `[l]  `^dl[dlll]  `^dr/14pt[dll]   [dll] \\ 
& {#5} \ar[r]  & {#6} \ar[r] & {#7}  \ar `r/10pt[d] `[l]  `^dl[dlll]  `^dr/14pt[dll]   [dll] \\
&  {#8} \ar[r] & ... & & }}
\newcommand{\rs}[5]{
\xymatrix{
         &            &            &  0         & \\ 
         &            &            & {#5} \ar[u] & \\
0 \ar[r] & {#1} \ar[r] & {#2} \ar[r] & {#3} \ar[r] \ar[u] & 0 \\  
         &            &            & {#4} \ar[u]\ar@{.>}[ul] & \\
         &            &            &  0  \ar[u] &
}}

\newcommand{\ts}[4]
{\xymatrix{
0 \ar[r] & {#1} \ar[r] & {#2} \ar[r] & {#3} \ar[r] & 0 \\  
         &            & {#4} \ar[u] \ar@{.>}[ul]^{\gamma} \ar@{.>}[ur]_{\lambda}  &            & \\
         &            &  0  \ar[u] &            &
}}

\newtheorem{theo}{Theorem}[section]

\newtheorem{prop}[theo]{Proposition}

\newtheorem{conj}{Conjecture}[section]

{\theorembodyfont{\rmfamily} }
{\theorembodyfont{\rmfamily} }

\numberwithin{equation}{section}

\begin{document}

\begin{titlepage}

\vspace{-5cm}

\title{
   \hfill{\normalsize UPR-1017-T} \\[1em] {\LARGE 
Branes, Bundles and Attractors:\\
Bogomolov and Beyond}
\\
[1em] }
\author{
     Michael R. Douglas$^{1,\&}$,
     Ren\'e Reinbacher$^1$
     and Shing-Tung Yau$^2$ \\[0.5em] {\normalsize
     $^1$Department of Physics and Astronomy, Rutgers University}
     \\[-0.4em] {\normalsize Piscataway, NJ 08854}\\[0.5em] {\normalsize
     $^2$Department of Mathematics, Harvard University}
     \\[-0.4em] {\normalsize Boston, MA 02138} \\[0.5em]
     {\normalsize $^\&$ I.H.E.S.,
        Le Bois-Marie, Bures-sur-Yvette, 91440 France}}

\date{}
\maketitle
\begin{abstract}
\noindent
We discuss conjectures following from the attractor mechanism in type
$II$ string theory about the possible Chern classes of stable holomorphic
vector bundles on Calabi-Yau threefolds.  In particular, we give  sufficient conditions for Chern classes to correspond to stable bundles.
\end{abstract}

\thispagestyle{empty}

\end{titlepage}

\section{Introduction} 

Let $(M,g)$ be a compact complex $n$-dimensional Kaehler manifold,
with Kaehler class $J\equiv [g]\in H^{(1,1)}_{\bar{\partial}}(M)$.
Let $V$ be a holomorphic vector bundle on $M$.  
The topological classification of such $V$ is given by a K theory class,
which up to torsion is determined by the Chern character $\ch(V)$.  
We recall that $\ch(V)$ is a sum of differential
forms in $H^{2*}(M)\equiv \oplus_k H^{2k}(M,\bbr)$, which given a connection 
with curvature $F\in H^2(M,\End V)$ can be written as
$$
\nch(V) = \Tr e^F
$$
Thus, its zero-form component is the rank of $V$, its two-form
component is the first Chern class $c_1$, and higher forms are
polynomial in the Chern classes.  We further recall that the bundle
$V$ is $\mu$-stable if, for all holomorphic subbundles $L$, we have
\begin{equation}
\label{eq-mustable}
\mu(L) < \mu(V) ,
\end{equation}
where the slope $\mu(E)$ of a bundle $E$ is defined by
\begin{equation}
\label{eq-slope}
\mu(E)=\frac{c_1(E)\cdot J^{n-1}}{r} .
\end{equation}
Note that a bundle is polystable if it is a direct sum of stable bundles with the same slope.
The same definitions can be made for a coherent sheaf $V$.
In addition, given a submanifold $D\subset M$ with embedding
$i:D\rightarrow M$, we  
can similarly define a $\mu$-stable holomorphic bundle $W$ on $D$.
In this case, we define $\nch(W)$ to be the class $i_{*}\Tr e^F$ on $M$ of
the coherent sheaf $i_{*}W$.

Now, given $J$, we define the set $\CH(J)\subset H^{2*}(M,\bbz)$ as the
subset of Chern characters for which there exists a $\mu$- polystable sheaf
$V$ or $W$.
We also define $\CHu \equiv \bigcup_J \CH(J)$, the subset of Chern
characters which support $\mu$- polystable sheaves for some $J$, and $\CHi
\equiv \bigcap_J \CH(J)$, the subset of Chern characters which support
$\mu$- polystable sheaves for all $J$.

These sets are known for curves and for a few surfaces.  Even for the
simplest surfaces, say $M=\cp2$, their precise structure is rather
intricate \cite{LePotier}.  Thus an approximate description, say by
upper and lower bounds, is also of interest.

In this work we discuss the problem of characterizing these sets, and
give some conjectures inspired by superstring theory.
Let us state our main conjectures without further ado:
\begin{conj}{
Let $M$ be a simply connected complex threefold $M$ with trivial
canonical bundle.  Consider
four classes
\begin{equation}
\label{conj-exCY}
r \in \mathbb{Z},\;r >1,\;\;\;c_i \in H^{2i}(X,\mathbb{Z}),\;\;\;i=1,2,3
\end{equation}
{and an ample  class  $\tilde{H}$  such that}
\begin{eqnarray}
\label{cond1}
\frac{1}{2r^2}\left(2r c_2 - \left(r-1\right)c_1^2-\frac{r^2}{12}c_2(X)\right)&=&\tilde{H}^2\label{condit2}\\
\frac{1}{6r^2}\left( c_1^3+3r(rch_3-ch_2c_1) \right)&< &\frac{2^{5/2}}{3} r\cdot \tH^3,\\
\end{eqnarray}
then, there exists a stable reflexive sheaf $V$ with respect to some ample class such that 
\begin{equation}
\rk(V)=r,\;\;\;c_i(V)=c_i\;\;\;i=1,2,3.
\end{equation}
In particular, there exists such a stable sheaf for any given ample class.}
\end{conj}

\begin{conj}{
Consider an smooth ample divisor $D$ on a simply connected Calabi-Yau threefold $X$  and three classes
\begin{equation}
\label{conj-exSur}
r \in \mathbb{Z},\;r>1,\;\;\;c_i \in H^{2i}(D,\mathbb{Z}),\;\;\;i=1,2
\end{equation}
{where $c_1$ is in the image of $H^{(1,1)}(X)\to H^{(1,1)}(D)$ such that}
\begin{equation}
\left( 2rc_2-\left(r-1\right)c_1^2-\frac{r^2}{12}c_2(D)\right)>0.
\end{equation}
Then there exists a stable vector bundle $V$ on $D$ such that 
\begin{equation}
\rk(V)=r,\;\;\;c_i(V)=c_i\;\;\;i=1,2
\end{equation}
In particular, there exists such a stable bundle for any given ample class.}
\end{conj}
The reason that stable reflexive sheaves arise here, along with a
more detailed discussion of stable objects in string theory,
will be given in section~\ref{sectionosm}.

We will give a ``physics proof'' of these conjectures in section 4.
Such conjectures arise in the study of compactification
of superstring theory, in more than one way.  The original problem 
of this type was that of
constructing realistic vacua from the heterotic string
\cite{CHSW};
see  \cite{DOPW,DOPR,Bouchard:2005ag}
and references there for the most recent work.
Here, one seeks a holomorphic bundle $V$ carrying a connection which
solves the Hermitian Yang-Mills equations, 
\begin{equation}
\label{eq-hermiteYM 2}
g^{a\bar{b}}F_{a\bar{b}}=c\mathbb{I} ,
\end{equation}
and satisfying the following topological constraints:
the rank is 3, 4 or 5; the first Chern class is vanishing, and the
second Chern class is determined by anomaly cancellation to be equal
to that of the tangent bundle of $M$.  
By the DUY (Donaldson-Uhlenbeck-Yau) theorem,
the question of whether such a bundle exists, is equivalent
to asking whether $\CHu$ contains an element satisfying these constraints.

Note that the topological conditions determine all but the third Chern
class $c_3$, which is unconstrained {\it a priori}.  Thus, knowledge
of $\CH$ would give us a list of allowed $c_3$ values.  Physically,
such a value determines the number of families of quarks and leptons
in four dimensions. Hence, a bound on the third Chern class in terms
of the other data translates into a interesting bound on the number of
families one can obtain in such compactifications.

A second physical problem, which will be of more primary interest in
this paper, is that of characterizing the BPS particles in type IIa
superstring theory, which arise from wrapping Dirichlet six-branes on
$M$ \cite{P}.  Loosely speaking, a Dirichlet
six-brane is a seven real dimensional minimal volume submanifold of
$M\times \bbr^{3,1}$, where $\bbr^{3,1}$ is Minkowski space-time,
carrying a vector bundle $V$ with a connection satisfying the
Yang-Mills equations.  A BPS particle is the special case in which the
embedding is a direct product of $M$ with a timelike geodesic in
$\bbr^{3,1}$, and the connection satisfies the Hermitian-Yang-Mills
equations.  The Chern character $\ch(V)$ determines the electric
and magnetic charges of this particle, and thus we find a second physical
interpretation of the set $\CHu$ -- it is the set of possible charges of BPS
particles.

The definition of $\mu$-stability, while precise, is not so easy to check,
as we must know all holomorphic subbundles of $V$.
Other necessary conditions for stability are known.  Perhaps the most
famous is the Bogomolov bound, which states that the Chern classes of a
$\mu$-stable holomorphic vector bundle $V$ will satisfy
\begin{equation}
\label{eq-Bogomolov}
(\frac{2c_2(V)}{r}-(r-1)\left(\frac{c_1}{r}\right)^2(V))
\cdot J^{n-2}\geqslant 0.
\end{equation} 
This follows from the existence of a solution to hermitian Yang-Mills,
and the positivity of the volume form.  

For surfaces, the Bogomolov bound already gives a rough picture of
$\CH(J)$.  Consider the example of the projective plane $\mathbb{P}^2$;
the actual boundary of $\CH(J)$ as given in \cite{LePotier2}, is given
by an inequality of the form \eq{Bogomolov} with a finite correction
depending on $\mu$.

Not much seems to be known for higher dimensions.  It is known that
any subset of $\CH(J)$ with fixed rank, $c_1$ and $c_2$ is finite
\cite{Ma,La}.

How can the connections we discussed to string theory help us with
this problem?  Our approach will be based on the attractor mechanism
\cite{FKS,FGK,Mor}.  This uses the relation between BPS particles and
extremal black hole solutions of supergravity theory.  The definition
of a BPS particle in type II supergravity provides a map $Z$ from the
Chern character and Kaehler class to $\bbc$, called the ``central
charge,'' \cite{Asp} explicitly
\begin{equation}
\label{eq-centralchargeBbrane}
Z[{\ch(V)}]= \int_M{\exp^{-(B+iJ)}}\ch(V) \sqrt{\td(X)}.
\end{equation}
Note that the Kaehler form $J$ has been promoted to
$B+iJ$, the ``complexified K\"ahler form'' in $H^2(M,\bbc)$.
%
%

Now, for a fixed $\ch(V)$, one can consider $Z$ as a function of
$B+iJ$, and look for a local minimum of the quantity
\begin{equation}
\label{eq-mass}
||Z||^2 \equiv \frac{|Z|^2}{\int_M J^3} .
\end{equation}
By definition, the minimum $||Z_{min}||$ is attained at
an ``attractor point'' for $V$ in complexified
Kaehler moduli space.  
We refer to a $\ch(V)$ for which $||Z_{min}||>0$ as satisfying the
``attractor condition,'' and let $\ATT\subset H^{2*}(M)$ be the set of 
points satisfying this condition.

In \cite{Mor}, Moore conjectured that if a charge vector is in $\ATT$,
then there will exist a stable BPS particle with that charge.
Specializing this conjecture to the case at hand leads to our basic
claim: that $\ATT\subset\CHu$.  Determining $\ATT$ requires no
knowledge about bundles and subbundles, only the de Rham
cohomology (with its product structure) and K\"ahler cone of $M$.

Working this out leads to a set of conditions
on Chern characters which are sufficient for existence.
The first of these
will turn out to be Bogomolov bound \eq{Bogomolov}.  
The next is a condition on the third Chern class:
if $V$ is stable and has vanishing first Chern class, and
$$
c_2(V)-\frac{\rk(V)}{24}c_2(M)=\tilde{H}^2
$$
for some ample Kaehler class $\tilde{H}$, we find
\begin{equation}
\label{eq-cthreebound}
|c_3(V)|\leqslant {\rm constant}\rk(V)\cdot{\tilde{H}^3},
\end{equation}
with an order one constant whose precise definition 
depends on stringy corrections; neglecting these it is $2^{5/2}/3$.

Let us now describe the attractor arguments in an intuitive way,
referring to \cite{Mor, De} for more details.  Given a Chern
character which satisfies the attractor condition, we want to see why
a BPS particle with these charges must exist.  This is because one can
find a extremal black hole solution of supergravity with these
charges, which to a distant observer is indistinguishable from a BPS
particle.  A more precise version of this argument relies on variation
of parameters (the dilaton field) to interpolate between the black
hole and the Dirichlet brane whose existence we are trying to
determine.  One can show that the existence of a BPS particle is
independent of this parameter (since the dilaton sits in a
hypermultiplet) and thus if the black hole solution is physically
sensible, so must be the Dirichlet brane configuration, implying that
the bundle and hermitian Yang-Mills connection of interest must exist.

Conversely, if $||Z_{min}||=0$, one finds a singular solution of
supergravity.  Here there are two subcases.  If the minimum is
achieved at a regular point in moduli space, the corresponding 
supergravity solution is clearly unphysical, and a BPS state of this
(single centered) type cannot exist.  On the other hand, if it is
achieved at a singular point (more precisely a boundary) of 
moduli space, supergravity breaks down, and one gets no clear statement
(in known examples, there is a BPS state).

The upshot is that the attractor condition is a sufficient but not
necessary condition for $\ch(V)\in \CH(J)$, for all $J$ in the basin
of attraction of the attractor point, and thus for $\ch(V)\in \CHu$.

To get a candidate for a necessary condition, and better understand
variation of stability, one must continue and discuss the far-reaching
generalizations of the original attractor conjecture made by Denef in
\cite{De,Denef:2001xn,Denef:2001ix,Denef:2002ru,Bates:2003vx}.  The
original attractor argument assumed that the supergravity solution was
spherically symmetric.  This need not be the case; more complicated
multi-center solutions provide examples of BPS bound states with
charges which do not satisfy the the original attractor condition.  In
general, the existence of these solutions is determined by the
existence of ``split attractor'' flows, this condition is manifestly
dependent on $J$ and has been shown to describe variation of stability
in examples.

To the extent that variation of stability is always described by this
mechanism, we would conclude that $\ATT\subset\CHi$.  
It also suggests a prescription for how to enlarge $\ATT$ to better
approximate $\CH$, as we discuss in section 5.

In order to trust these physical arguments, we need two conditions
to be satisfied.
One of these, to justify the use of the large volume expression for
the central charge \eq{centralchargeBbrane} and the comparison with
$\mu$-stability, is that the attractor point satisfy $J >> 1$ in string
units (more precisely, that world-sheet instanton corrections be
negligible).  We will discuss this condition later on; in principle the
discussion could be generalized to include world-sheet instantons and
the more general notions of stability introduced in \cite{Do,Br}.

The other condition, required so that supergravity is a good
approximation to the more exact string theory description, is that the
quantity $||Z_{min}||^2 >> 1$, as this is the condition that the area
of the event horizon (and thus the curvature radius) be large compared
to the Planck scale.  A sufficient condition for this is for all
non-zero components of the Chern character $\ch(V)$ to be large.  

More precisely, since \eq{mass} is homogeneous of degree two in
$\ch(V)$, for any given $\ch(V)$, there will exist an $N_{min}$ such
that for all $N>N_{min}$ the set of Chern characters obtained by the
rescaling $\ch(V) \rightarrow N \ch(V)$ will satisfy this condition.
Note that the definition of attractor point and the attractor
condition are unchanged by such an overall rescaling.

We will refer to subsets of $H^{2*}(M,\bbz)$ or bounds on $\ch(V)$ which
are invariant under rescaling as ``homogeneous.''  To the extent that
the set of charges of BPS objects and therefore the sets $\CH$ can be
understood from the genus zero attractor considerations we discuss, we
conclude that these sets are homogeneous.

This condition plays little role in the existing mathematical
discussions, and it is not immediately obvious that non-trivial
homogeneous bounds exist.  Note that it is important that we are asking
for invariance under rescaling of the Chern character and not the
Chern classes; for example the Bogomolov bound \eq{Bogomolov} is 
not invariant under rescaling the $c_i$ and $r$, but upon rewriting
it as
\be\label{eq-bbound}
\Delta_2 \equiv
 \left(\frac{1}{2}\left(\frac{c_1}{r}\right)^2
 - \frac{ch_2}{r}\right)\cdot J^{n-2} \ge 0 ,
\ee
we see that it is homogeneous.  Thus we raise the possibility that,
along with our sufficient conditions, homogeneous necessary conditions
on the Chern characters also exist.  We discuss this in section 6.

Before making this physical discussion a bit more precise and work out
the attractor bounds, we will in the next section, state our results
and the related mathematical conjectures on the existence of stable
bundles on Calabi-Yau threefolds and surfaces with ample canonical
bundle and test them against hard mathematical theorems.

\noindent
\vskip 0.2in
Note that the previous version of the paper contained a possible conjecture on a slight strengthening of Bogomolov inequality. However, as pointed out by M. Jardim,  there exists examples which violate the stronger inequality. We will present an example due to M. Jardim in Appendix  $B$.

\section{Existing mathematical results}

We begin this section by setting out some elementary features of the
problem, and then briefly summarize some known mathematical existence
theorems and constructions.

Let $\cH\cong H^{2*}(M,\bbr)$ be the direct sum of the even de Rham
cohomology groups.  We consider it as a commutative algebra, with the
standard cohomology product, and with the fundamental class
$\int_M:\cH\rightarrow\bbr$.  Should we need it, a basis for
$H^2(M,\bbr)$ will be denoted $\{\omega^i\}$ with $1\le i\le b^2(M)$.
We can take these to be ample classes if we like.

Of course, only classes in $H^{p,p}(M,\bbc)$ can be realized as
the Chern characters of holomorphic bundles.  Our primary interest is
in $M$ a Calabi-Yau threefold with $h^{1,0}=h^{2,0}=0$, for which these
span $\cH$.

The Mori cone $MC\in H_2(M,\bbz)$ is generated by effective curves
in $M$.
Its dual  $KC\in H^2(M,\bbr)$ of ample classes $\omega$ satisfying
$<\omega,C> \;\ge 0 \;\forall C\in MC$ is the K\"ahler cone.

Ample classes have various positivity properties \cite{Lar}.  We will
use the following: consider $H_i$ with $i=1,2,3$ ample on an
irreducible complete threefold (such as $M$); then
\be
\label{eq-ample}
(H_1H_2H_3)^3 \ge H_1^3 H_2^3 H_3^3.
\ee

Given a sheaf $V$ on $M$, we define the generalized Mukai vector
$\gamma(V)\in \cH$ as
\be \label{eq-mukai}
\gamma(V) = \Tr e^F \sqrt{{\rm Td~ M}} .
\ee
This is of course simply related to $\nch(V)$ by a linear change of
basis, but turns out to emerge more naturally from the physics.
One of its mathematical advantages is a simple statement of the
Grothendieck-Riemann-Roch formula.  Defining an involution on $\cH$ 
which corresponds to dualizing a vector bundle,
$$
(-)^k : \omega^{2k} \rightarrow (-1)^k \omega^{2k} ,
$$
we have
\be \label{eq-intersect}
\chi(V,W) = \sum (-)^l \dim {\rm Ext}^l(V,W) =
 \int_M \gamma(V) \wedge (-)^k \gamma(W) .
\ee
As this depends only on the Mukai vectors $\gamma(V)$ and $\gamma(W)$,
we also write this as $\chi(\gamma(V),\gamma(W))$.

We now consider the subsets $\CH$, $\CHi$ and $\CHu$ of the introduction,
defined using the map \eq{mukai}.  All of these contain the Mukai vector
of the trivial bundle, $\gamma(\om)=\sqrt{\td M}$, as well as of all other
line bundles.  More generally, since $\mu$-stability is invariant
under tensoring by a line bundle, we have
$$
e^F ~ \CH \cong \CH \qquad \forall F\in H^2(M,\bbz)
$$
(resp. $\CHi$ and $\CHu$).

Indeed, since the tensor product of two $\mu$-polystable bundles
(at a fixed $J$) is $\mu$-polystable, each $\CH(J)$ is closed under
multiplication -- it is a subalgebra of $\cH$.

\subsection{Homogeneous invariants} 

One way to simplify the problem is to discuss asymptotic results, valid
for large Chern character.  Thus, we define $\PCH$,
$\PCHi$ and $\PCHu$ to be the closure of the projectivized versions of
the above; in other words $v\in\PCH$ if for any $\delta$
we have $\lambda v+\epsilon\in\CH$ for
some $\lambda\in\bbr^+$ and some $\epsilon\in\cH$ with 
$|\epsilon|<\delta$, etc.

One can define polynomials in the Chern classes which are
invariant under an overall rescaling of $\ch(V)$, and are thus
well defined on $\PCH$.   A useful set of generators for these, introduced by
Drezet, is defined by the following expansion,
\begin{equation} \label{eq-drezet}
\log \ch(V) = \log r + \sum_{i=1}^n (-1)^{i+1} \Delta_i(V),
\end{equation}
where $\Delta_i(V)$ is a $2i$-form.
For example, $\Delta_1=c_1/r$, $\Delta_2$ is as in \eq{bbound}, 
and $\Delta_3 = c_3/2r + c_1(\ldots)$.  
The $\Delta_i$ for $i\ge 2$ 
are unchanged by tensoring with line bundles, and
this can be used to infer the terms proportional to $c_1$.

We furthermore have
$$
\Delta_i(E\otimes F) = \Delta_i(E) + \Delta_i(F)\; \forall\; 1\le i\le n ,
$$
and in this sense the $\Delta_i$ linearize the algebra structure on $\cH$.
Thus we have
\begin{prop}
Let $\Delta\CH(J)$ be the image of $\CH(J)$ in $\cH /(H^0(M)\oplus H^2(M))$  under the map 
$\Delta: \cH \rightarrow (\Delta_2,\ldots,\Delta_n)$; 
it is a convex cone, invariant under $\Delta_i\rightarrow (-1)^i\Delta_i$. In particular, note that $\Delta$ factorizes through
$$
\xymatrix{
\CH(J) \ar[r]^-{\Delta}\ar[d] &\cH /(H^0(M)\oplus H^2(M))\\
\PCH(J) \ar[ru]^-{\Delta}
}
$$
\end{prop}

For example, the Bogomolov bound is
\be\label{eq-bboundtwo}
\Delta_2 \equiv
 \left(\frac{1}{2}\left(\frac{c_1}{r}\right)^2
 - \frac{ch_2}{r}\right)\cdot J^{n-2} \ge 0 ,
\ee
which is consistent with this.
For $n=3$, the set $\Delta\CH(J)$ is determined by one more inequality,
$$
|\Delta_3| \le f(\Delta_2)
$$
for some convex function $f$.  For example, 
we might have $f(x)=x^a$ for any $a\ge 1$.

\subsection{Existence results}
\label{ss-suff}

We consider the moduli spaces $\mathcal{M_{J}}(r,c_1,c_2)$
and $\tilde{\mathcal{M}}_{J}(r,c_1,c_2) $ of $\mu$-stable vector
bundles and sheaves respectively (with respect to $J$), with rank $r$
and first and second Chern classes $c_1$ and $c_2$.  If 
$\mathcal{M_{J}}(r,c_1,c_2)$ is not empty, we say that bundles of
(resp. sheaves) of this topological type exist.

First, for a general smooth algebraic variety $M$ of arbitrary dimension,
an asymptotic result of Maruyama \cite{Ma} states that bundles exist
for $r\geqslant \dim M\geqslant2$ and $ c_2(V) \omega^{n-2} \gg 0$.

More precise statements are available for algebraic
surfaces. An early general result (Dr\'ezet-Le Potier) is that
sheaves exist on $\mathbb{P}_2$ if
\begin{equation}
2rc_2-(r-1)c_1^2 \geqslant  2r^2\delta,
\end{equation}
where $\delta$ is a periodic function $1/2 \leqslant \delta
\leqslant1$. Actually, this result is sharp, that is, the imposed
condition is also necessary if one excludes exceptional bundles, that
is bundles with $Ext^1(W,W)=0$.

On a $K3$ surface, we have the following recent result
by Yoshioka \cite{Yos}: semistable sheaves (for  a general ample class) exist if
\begin{equation}
\label{eq-Yos}
2rc_2-(r-1)c_1^2-\frac{r^2}{12}c_2( K3) \geqslant -2.
\end{equation}
This is simply the requirement that the Mukai vector 
$v(W)=\ch(W)\sqrt{\td(K3)}$ satisfies $(v,v)\geqslant-2$. Under the assumption that $v$ is primitive, Yoshioka can actually show the existence of stable sheaves. After the exclusion of the rank one case, these sheaves are generically locally free.

For general surfaces there are existence results by Taubes
\cite{Taubes}, Gieseker \cite{Gieseker}, Artamkin \cite{Artamkin},
Friedmann, J. Li, etc. In particular, bundles with $r=2$ and $c_1$
exist for
\begin{equation}
c_2 \geqslant 2h^2(D,\mathcal{O}_D)+2.  
\end{equation}
More generally, W. P. Li and Z.Qin have shown that
$\mu$-stable bundles with given $r$ and $c_1$ exist if
\begin{equation}
c_2>\alpha,
\end{equation}
where $\alpha$ is a numerical constant depending on the surface $D$,
$r,\;c_1$ and the class $J$. Also, under the assumption $c_1\neq 0$, constant $\alpha$ can not be universal depending on $r$ and $c_1$ only. It must, in general, depend on $J$. In particular, on Hirzebruch surfaces, there exists a $c_1$ such that for any $c_2$ we can find an ample class $J$ such that  $\mathcal{M_{J}}(r,c_1,c_2)$ is empty.

A stable bundle $W$ on a hypersurface $D$ in $M$ gives us a point in $\CH$,
by using the Grothendieck-Riemann-Roch formula:
\begin{equation}
\ch(i_{*}W)=i_{*}(\ch(W)\td(D))\td(M)^{-1},
\end{equation}
where $i: D \subset M$ denotes the inclusion. Using
\begin{equation}
c_1(D)=-D^2,\;\;\;c_2(D)=D^3+c_2(M)D
\end{equation}
we find
\begin{equation}
\begin{split}
\rk(i_{*}W)&=0\\
ch_1(i_{*}W)&=\rk(W)D\\
ch_2(i_{*}W)&=-\frac{\rk(W)D^2}{2}+c_1(W)\\
ch_3(i_{*}W)&=\frac{\rk(W)D^3}{6}+ch_2(W)-\frac{(D^2) \cdot  c_1(W)}{2}.\\
\end{split}
\end{equation}
Thus, the general import of the results above is that $\PCH(J)$ and
$\PCHu$ contain the points $1-e^{-D}-c_3 J^3$, for any effective
divisor $D$ and for sufficiently large $c_3$.

No such general existence theorems for stable vector bundles on
Calabi-Yau threefolds are known.  Of course, we know the tangent and cotangent
bundles are stable, as are tensor products of these.  
Two other well known constructions are the monad construction and the spectral cover
construction.

As an example, the Chern character for $TQ$, the tangent bundle of the quintic
hypersurface, is
\begin{equation}
\label{eq-chTQ}
\ch(TQ)=3-10H^2_{Q}-20H^3_Q,
\end{equation}
where we denote the restriction of the hyperplane class to the quintic by $H_Q$. 

The monad construction was used by Maruyama in \cite{Ma} to prove his
asymptotic result. We follow here the presentation of \cite{Dis}.  To
begin with, fix an integral Kaehler class $H$ and a smooth hypersurface
$D$ in $X$, with $[D] \cong r[H]$ in cohomology.
Consider the exact sequence on $D$,
\begin{equation}
\xymatrix{
{\mathcal{O}}_D(-nH) \ar[r]^-{s_i}& {\mathcal{O}}_D^{\oplus r}
 \ar[r] ^{\alpha}&{\mathcal{F}}},
\end{equation}
where $n$ and $r$ are some positive integer ($r\geqslant 3$) and
${s_i}$ are $r$ linearly independent global sections of
$\mathcal{O}_D(nH)$ which are base point free. This implies that
$\mathcal{F}$ is a vector bundle on $D$. Composing the restriction map
$r$ from $M$ to $D$ with $\alpha$ gives the vector bundle
$\tilde{E}(-D)$ on $X$ as
\begin{equation}
\xymatrix{
\tilde{E}(-D) \ar[r]& {\mathcal{O}}_X^{\oplus r} \ar[r]
 ^{\alpha\circ r}&{\mathcal{F}}}.
\end{equation}
Provided that $n$ is sufficiently large, it was shown in \cite{Ma}
that $\tilde{E}(-D)$ is stable. 
We are interested in the stable rank
$r$ vector bundle $E=\tilde{E}(-J)$. Its Chern classes were computed
in \cite{Dis}.
\begin{equation}
\label{eq-monadclass}
\begin{split}
c_1(E)&=0\\
c_2(E)&=\frac{r}{2}(2n+1-r)H^2\\
c_3(E)&=\left(-\frac{1}{6}r(r-1)(r-2)+
 \frac{r}{2}\left(2(n+1)^2-r(2n-r+3)\right)\right)H^3
\end{split}
\end{equation}

The spectral cover construction \cite{don,fmw} applies only for
elliptically fibered varieties.  We consider an
elliptically fibered Calabi-Yau threefold $\pi: M \to B$, and
construct a rank $r$ bundle $V$ with $c_1=0$,
stable respect to some ample class of the form
\begin{equation}
H=H_0+n\pi^{*}\alpha,\;n \gg0,\;H_0,\alpha\; {\rm ample}.
\end{equation}
The second Chern class of the the tangent bundle of these threefolds
and the second Chern class of these bundles are given by
\begin{equation}
\begin{split}
c_2(M)&=12\sigma \pi^{*}c_1(B)+m_M F\\
c_2(V)&=\sigma \pi^{*}{\eta}+m_V F\\
\end{split}
\end{equation}
where $\sigma$ and $F $ denote the zero section and fiber class
respectively, $\eta$ some ample class on $B$, and $m_M$ and $m_V$
are positive integers determined by the choice of threefold and the
choice of spectral cover. More specifically, the class of the spectral cover is
given by $[C]=r\sigma+\pi^{*}\eta$, with stability ensured for irreducible
spectral covers. The existence of such a cover can be guaranteed by
imposing the following numerical constraints
\begin{equation}
\eta-rc_1(B)\;\;\;\text{effective},\;\;\eta\;\;\;\text{ample}.
\end{equation}
These results give us a supply of known elements of $\CH$.  

\subsection{Necessary conditions}

We now consider necessary conditions on $\CH$.  First, we repeat
the Bogomolov bound for a $\mu$-stable bundle with respect to $J$,
\be
\Delta_2 \equiv
 \left(\frac{1}{2}\left(\frac{c_1}{r}\right)^2
 - \frac{ch_2}{r}\right)\cdot J^{n-2} \ge 0 .
\ee
For Calabi-Yau threefolds this seems to be sharp.  
On surfaces, evidence so far is consistent with a slightly stronger bound. We postpone the discussion to Appendix A which will fit well
with the physics arguments below.

Bounds on the higher Chern classes are discussed in \cite{La2}.
For example, given  any polynomial $P(c)$ in the Chern classes,
there exist computable polynomials $P_L$ and $P_H$ in $r$, $c_1$,
$\Delta$ and invariants of $M$, such that 
\be\label{eq-gbound}
P_L \le P(c) J^{n-k} \le P_H
\ee
for all semistable reflexive sheaves with fixed $r$, $c_1 J^{n-1}$
and $0\le c_2 J^{n-2}\le \Delta$.

In terms of the invariants \eq{drezet}, the general bound \eq{gbound} becomes
$$
|\Delta_i(V) J^{n-1}| \le W_i(\Delta_2(V)J^{n-2})
$$
in terms of computable polynomials $W_i$ (which again depend on invariants
of $M$).

\section{Mirror symmetry and the attractor mechanism}

In this section, which is not strictly necessary for the rest of the
discussion, we explain the relation between the IIa string theory in
which our problem is naturally formulated, and the mirror IIb string
theory.  This is useful in part because some of the definitions are
more natural in the IIb picture, and in part to explain the string
theory corrections.

Thus, we consider the IIb superstring theory \cite{GSW}, compactified
on $\tilde{ M}\times\bbr^{3,1}$.  Later we will take $\tilde M$ to be
the mirror of $M$, but for now the discussion is general.  Such a
compactification leads to a four dimensional $N=2$ supergravity theory
containing various fields, in particular $b^3/2$ abelian gauge fields
(thus, copies of Maxwell's theory).  These arise because the IIb
supergravity theory contains a four-form generalized gauge potential
$C^{(4)}$, satisfying the self-duality relation (up to non-linear
terms we can ignore),
$$
dC^{(4)} = *dC^{(4)} ,
$$
where $*$ is the ten-dimensional Hodge star.  Then, 
given a class $\Sigma_i\in H_3(M,\bbz)$, the integrals
$$
A_i \equiv = \int_{\Sigma_i} C^{(4)}
$$ define a $b_3$-dimensional linear space of one-form abelian
connections.  The self-duality relation then implies that only half of
these are physically independent; the other half are dual in the sense
of four-dimensional electric-magnetic duality.  Thus, the lattice of
conserved electric and magnetic charges carried by particles in four
dimensional space-time is
\begin{equation}
\Lambda=H^3(\tilde{ M},\mathbb{Z}).
\end{equation}
Consider a charge $\hat{\gamma}\in H^3(\tilde{ M},\mathbb{Z})$, and
assume there exists a corresponding BPS state.
A natural origin for such a BPS state in IIB string theory is a Dirichlet 
3-brane.  It was shown
in \cite{BBS} that BPS D3-branes can be represented in the large volume
limit of $\tilde{ M}$ by a special Lagrangian three-cycle $\gamma_L$.

As outlined above, we will attempt to decide about the existence of BPS
states in the IIB super-gravity approximation of IIB string theory.
As shown in \cite{FKS}, such BPS states can be obtained as dyonic
black hole solutions on $M_4$. More precisely, we consider static,
spherically symmetric four dimensional space-time configurations which
are asymptotically flat and carry a dyonic charge $\hat{\gamma} \in
H^3(\tilde{ M},\mathbb{Z})$. Under these assumptions, the equations for
the BPS states reduce to a dynamical system on $\mathbb{R}^{*}\times
\tilde{\mathcal{M}}$ involving $Z_L$ and $|Z_L|$ defined in
\ref{eq-period} and \ref{ eq-periodnorm}. $\mathbb{R}^{*}$ describes
the radial direction of the spherically symmetric space-time which
serves as the evolution parameter and $\tilde{\mathcal{M}}$ the
universal cover of the moduli space $\pi: \mathcal{\chi} \to
\mathcal{M}_{Complex}(\tilde{ M})$ of complex structures of the
Calabi-Yau threefold $\tilde{ M}$.  Let us consider the family of
marked Calabi-Yau threefolds $\pi: \tilde{\mathcal{\chi}} \to
\tilde{\mathcal{M}}(\tilde{ M})$ and the pull-back
$\tilde{\mathcal{L}}$ of the Hodge bundle
$\mathcal{L}=R\pi_{*}\omega_{\chi/\mathcal{M}}$.  Then we define for
each $\gamma \in H_3(\tilde{ M},\mathbb{Z})$, $\Omega \in
\tilde{\mathcal{L}}$ a function on the total space
$\tilde{\mathcal{L}}\to \tilde{\mathcal{M}}$
\begin{equation}\label{eq-period}
Z_L(\gamma,\Omega) =  \int_{\gamma} \Omega
\end{equation}
and
\begin{equation}
\label{ eq-periodnorm}
|Z_L(\gamma,\Omega)|^2=\frac{|\int_{\gamma} \Omega|^2}{i \int_{\tilde{ M}} \Omega \wedge \bar{\Omega}}
\end{equation}
Following \cite{Mor}, the existence of a BPS states in the
super-gravity approximation can be stated in terms of
$|Z_L(\gamma,\Omega)|$. There are three cases:
\begin{enumerate}
  \item$|Z_L(\gamma,\Omega)|^2$ has a nonvanishing local minimum. In
  this case we expect to have a BPS state in the theory.
  \item It can happen that $|Z_L(\gamma,\Omega)|^2$ has no stationary
  point in $\tilde{M}$. It might or might not vanish at the
  boundary. In this case the supergravity approximation breaks down,
  and we can not decide whether $\gamma$ supports an BPS state or not.
  \item For some vectors $\gamma$ it might happen that
  $Z_L(\gamma,\Omega)=0$ for some complex structure in the interior
  of $\tilde{\mathcal{M}}$. Naively, such charges do not support BPS
  states.  However, there might be split attractor solutions which nevertheless
  are BPS states.
\end{enumerate}

The next step is to use the key result of \cite{FGK} which states that
$|Z_L(\gamma,\Omega)|^2$ has a stationary point at $z_{*}(\gamma) \in
\tilde{\mathcal{M}}$ (which is necessarily a local minimum), with
fixed point value $Z_{*}\neq 0$ iff Poincare dual of a given charge
$\gamma \in H_3(\tilde{ M},{\mathbb{Z}})$ can be written as
\begin{equation}\label{eq-attractor}
\hat{\gamma}=\hat{\gamma}^{(3,0)}+ \hat{\gamma}^{(0,3)} \in H^{(3,0)}(\tilde{ M}) \oplus H^{(0,3)}(\tilde{ M})
\end{equation}
Hence, if we can solve for the complex structure in the interior of
$\tilde{\mathcal{M}}$ such that \ref{eq-attractor} holds, the
existence of a BPS state of charge $\hat{\gamma}$ is shown. Equation
\ref{eq-attractor} are the celebrated attractor equations which were
first solved in large volume limit by \cite{Shmakova}

\subsection{Review of closed string mirror symmetry}\label{s-mirror}

We generally follow the reviews \cite{Asp,HKT}.  The strongest version of
mirror symmetry for the closed string states that, if $ M$ and
$\tilde{ M}$ are Calabi-Yau threefolds which form a mirror pair, IIA
string theory compactified on $ M$ is isomorphic to IIB string theory
compactified on $\tilde{ M}$.  This is believed to follow from the
corresponding statement about topological string theory, namely that
$ M$ and $\tilde{ M}$ are mirror pairs if the operator algebra of the
A-model with target space $\tilde{ M}$ is isomorphic to the operator
algebra of the B-model with target space $ M$.

To define the A(B) model on $\tilde{ M}( M)$, one considers a
two-dimensional non-linear sigma model with $(2,2)$ super-symmetry on
$\tilde{ M}( M)$. The A(B)-model can be obtained by twisting the
original superconformal field theory.

One finds that  the space of operators for the A-model is given by $H^{*}(\tilde{ M},\mathbb{C})$. Also, the A-model depends only on the complexified Kaehler form $(B+iJ)(\tilde{ M}) \in \mathcal{M}_{Kaehler}(\tilde{ M})$, where
\begin{equation}
\mathcal{M}_{Kaehler}(\tilde{ M})=(H^2(\tilde{ M},\mathbb{R}) + i\mathcal{K}_{\tilde{ M}})/ H^2(\tilde{ M},\mathbb{Z})
\end{equation}
and $\mathcal{K}_{\tilde{ M}}$ denotes the Kaehler cone of $\tilde{ M}$. More precisely, all deformations of the A-model are described by $(B+iJ)(\tilde{ M})$. The operator product in the large volume limit $(J \to \infty)$ agrees with the cohomology ring given by the wedge product on $H^{*}(\tilde{ M},\mathbb{C})$. 
For a generic point in the moduli space $\mathcal{M}_{Kaehler}(\tilde{ M})$ the operator  product will receive corrections by world sheet instantons, that is, by  holomorphic maps from $S^2 \to \tilde{ M}$, and is called quantum cohomology ring.

For the B-model on $ M$, one finds that the space of operators is given by the $\bar{\partial}$-cohomology $H^{*}_{\bar{\partial}}( M)$. Also, the B-model depends only on the complex structure $\mathcal{M}_{Complex}( M)$ of $ M$. More precisely, the deformations of the B-model are described by elements of $H^1_{\bar{\partial}}( M,TM)$.  The 3-point functions for elements of $H^1_{\bar{\partial}}( M,TM)$ is given by the canonical paring
\begin{equation}
H^1_{\bar{\partial}}( M,TM)\otimes H^1_{\bar{\partial}}( M,TM) \otimes H^1_{\bar{\partial}}( M,TM) \to H^3( M,\mathcal{K}_M) \cong H^0( M,\mathcal{O}_ M),
\end{equation}
and has no instanton corrections.

One of the first conclusions one can draw from mirror symmetry is on cohomological level, namely, that
\begin{equation}
h^{p,q}(\tilde{ M}) = h^{3-p,q}( M).
\end{equation}
In addition, since the A-model on $\tilde{ M}$ depends on the
complexified Kaehler moduli space $\mathcal{M}_{Kaehler}(\tilde{ M})$
and the B-model on $ M$ depends on the complex structure
$\mathcal{M}_{Complex}( M)$, a precise statement of mirror symmetry
must include a map between these two moduli spaces. This map is called
mirror map. However, this map is not globally defined. One can define
the map locally around a base point which is usually given by the
large radius limit $\tilde{ M}$ and map it to the large complex
structure limit of $ M$\footnote{This point might not be unique.}. In
particular, one finds on both sides so-called special coordinates.

For the B-model which depends on the complex structure $\mathcal{M}_{Complex}( M)$ of $ M$, these coordinates are easily described.  We choose a symplectic basis $
(\hat{\alpha}_I,\hat{\beta^I}) \in H^3( M,\mathbb{Z})$
such that
\begin{equation}\label{eq-symlectic}
\int_{ M} \hat{\alpha_J}\wedge \hat{\beta}^J=<\hat{\alpha}_J,\hat{\beta^I}>= \delta^I_J,\;\;\;<\hat{\alpha}_J,\hat{\alpha}_I>=0,\;\;\; <\hat{\beta}_J,\hat{\beta}_I>=0
\end{equation}
and define the Poincare dual basis $
(\alpha^I,\beta_I) \in H_3({ M},\mathbb{Z})$
such that 
\begin{equation}\label{}
\int_{\alpha^I} \hat{\alpha_J} =  \int_{{ M} }\hat{\alpha_J}\wedge \hat{\beta}^J,\;\;\;
\int_{\beta_I} \hat{\beta^J} =  \int_{{ M} }\hat{\beta^J}\wedge \hat{\alpha}_J.
\end{equation}
Now we can introduce the period vectors 
\begin{equation}\label{eq-period MF}
 X^I=\int_{\alpha^I} \Omega,\;\;\; F_I = \int_{\beta_I} \Omega.
\end{equation}
The coordinates $ X^I,I=1,\ldots,h^{2,1}( M)+1$ are called projective coordinates on $\mathcal{M}( M)$. In addition, it can be shown that there exists a so-called pre-potential $F$, a homogenous function of degree two in $X^I$ such that $
F_I = \frac{\partial F(X)}{\partial X^I}$.
It follows from \ref{eq-period} that the the holomorphic three-form $\Omega$ can be written as
\begin{equation}\label{eq-3form}
\Omega = \sum X^I \hat{\alpha}_I - \sum F_I \hat{\beta}^I.
\end{equation}
The large complex structure limit is a point in $\mathcal{ M}$ that will be determined by a vanishing $3$-cycle which admits a  maximal unipotent monodromy transformation on $H^3( M,\mathbb{Z})$. After possible relabeling we call the vanishing $3$-cycle $\alpha_0$. This choice 
singles out the coordinate $X^0$ and allows the introduction of
special coordinates
\begin{equation}
\label{eq-specialcoordinates}
t^a=\frac{X^a}{X^0},\;\;\;a=1,\ldots,h^{(2,1)}({ M}).
\end{equation}
These coordinates allow the definition of an in-homogenous pre-potential $\mathcal{F}(t)=\frac{F(X)}{(X^0)^2}$ which can be explicitly computed \cite{HKT} and is given to leading order (without instanton corrections)
\begin{equation}
\label{eq-inhomogenous}
\mathcal{F}(t)=\frac{1}{3!}\sum D_{abc}t^at^bt^c+\frac{1}{2}\sum A_{ab}t^at^b+\sum \mathrm{B}_at^a + \cdots
\end{equation}
In particular, knowing the in-homogenous pre-potential, allows to express the derivatives of $F(X)$ as 
\begin{equation}\label{eq-derivativeF(X)}
{F_0} =X_0 (2\mathcal{F}(t)-
\sum t^i \partial_i \mathcal{F}(t)),\,\,\,
{F_i} =  X_0\partial_i \mathcal{F}(t),\,\,\,I=(0,i).
\end{equation} 

The special coordinates of the A-model on $\tilde{ M}$ can be introduced by
\begin{equation}
\label{eq-specialkaehlercoordinates}
B+iJ=\sum_a \tilde{t}^{'a} B_a + i\tilde{t}^{''a}J_{a},
\end{equation}
where $\{B_a\}$ is basis of  $H^2({\tilde{ M},\mathbb{R}})$ and $\{J_a\}$ a framing of $\tilde{ M}$. In particular, the set $\{\tilde{t}^{a}=\tilde{t}^{'a}+i\tilde{t}^{''a}: t^{'a}=0,t^{''a}\geqslant 0\}$ is contained in the closure of the Kaehler cone of $\tilde{ M}$.
With the correct choice of these basis the isomorphism between $\mathcal{M}_{Kaehler}(\tilde{ M})$ and $\mathcal{M}_{complex}( M)$ around the large complex structure limit is given by 
\begin{equation}
\label{eq-isomorphis(t) }
t^a=\tilde{t}^{a},\;\;\;a=1,\cdots,h^{2,1}( M)=h^{1,1}(\tilde{ M})
\end{equation}
Also, the coefficients in the in-homogenous pre-potential \ref{eq-inhomogenous} are expressed easiest in terms of  $\tilde{ M}$. One finds \cite{HKT}  $(-D_{abc})$ to be the intersection matrix of
\begin{equation}
H^2(\tilde{ M},\mathbb{Z}) \times H^2(\tilde{ M},\mathbb{Z}) \times H^2(\tilde{ M},\mathbb{Z}) \to H^6(\tilde{ M},\mathbb{Z}),
\end{equation}
$A_{ab}$ a symmetric integral matrix and $\mathrm{B}_a$
\begin{equation}
\mathrm{B}_a=- \frac{1}{24}\int_{\tilde{ M}} J_a \wedge c_2(\tilde{ M}).
\end{equation}
We will conclude our discussion on mirror symmetry of the closed string by the following remark. Instead of the A(B)-model on $\tilde{ M}( M)$ we can consider the A(B)-model on $ M(\tilde{ M})$. This will give an local identification of $\mathcal{M}_{Kaehler}( M)$ and $\mathcal{M}_{complex}(\tilde{ M})$ in the large complex structure limit of $\tilde{ M}$. This isomorphism is given by
\begin{equation}
\label{eq-isomorphism(s) }
s^{\alpha}=\tilde{s}^{\alpha},\;\;\;\alpha=1,\ldots,h^{2,1}(\tilde{ M})=h^{1,1}( M),
\end{equation}
where $\{s^{\alpha}\}$ are the special coordinates on
$\mathcal{M}_{complex}(\tilde{ M})$ introduced in
\ref{eq-specialcoordinates} and $\{\tilde{s}^{\alpha}\}$ the special
coordinates on $\mathcal{M}_{Kaehler}( M)$\footnote{We will always
denote the coordinates on the Kaehler moduli space with a tilde, and
the coordinates on the complex moduli space without a tilde.}
introduced in \ref{eq-specialkaehlercoordinates}.

\subsection{Dirichlet branes and open string mirror symmetry}
\label{sectionosm}

Our main point here will be to review how the BPS branes we have 
just discussed are realized as concrete objects in string theory,
special Lagrangian $3$-cycles on $\tilde M$ in the IIb theory, 
and holomorphic objects on $M$ in IIa theory,
where $M$ is the type II mirror of $\tilde{ M}$.
For reviews of this subject, consult \cite{Asp,HKT,Sharpe} and the upcoming
\cite{Clay}.  

To define open strings, we introduce Dirichlet branes, which in the
first instance are submanifolds on which an open string can end,
carrying vector bundles.  Starting with this definition, it can be
shown that the only D-branes compatible with the A-model are
Lagrangian 3-cycles, while D-branes compatible with the B-model are
holomorphic submanifolds.  A more detailed analysis shows, that
A-branes on $\tilde{ M}$ are described by the Fukaya category
$\mathcal{FU}({\tilde{ M}})$ and B-branes on $ M$ by the derived
category of coherent sheaves $\mathcal{D}( M)$. The mirror symmetry
conjecture by Kontsevich states that for a mirror pair $ M$ and
$\tilde{ M}$, the categories $\mathcal{D}( M)$ and
$\mathcal{FU}({\tilde{ M}})$ are equivalent\footnote{Note that there
are subtleties if $\mathcal{FU}({\tilde{ M}})$ is not triangulated.}.
We are interested in the A and B branes of $\tilde{ M}$ and $ M$ (in the
large volume limit of $\tilde{ M}$ and $ M$) which descent from
BPS-states in the untwisted theory.

To begin with, we consider A-branes on $\tilde{ M}$.  Before twisting, these are
special Lagrangian 3-cycles $\gamma_L$, for which the
holomorphic 3-form $\Omega$ obeys
\begin{equation}
\Omega|_{\gamma_L}= \mathrm{ constant}.
\end{equation}
The electric and magnetic charge of $\gamma_L$ is given by
its Poincare dual $\hat{\gamma}_L \in H^3(\tilde{ M},\mathbb{Z})$, and
its central charge by
\begin{equation}\label{eq-centralchargeIIB1}
\mathcal{Z}_{\gamma_L} = \int_{\gamma_L} \Omega.
\end{equation}
In particular, note that the central charge of the A-cycles depends on
the complex structure of $\tilde{ M}$, while A-model correlation functions depend only
on the Kaehler structure $(B+iJ)(\tilde{ M})$.

In general, a B-brane on $ M$ which descends from a BPS state in
the untwisted theory corresponds to a $\Pi$-stable object
$\mathcal{E}^{*}$ in the derived category $D(\Coh~ M)$.
Such an object is a quasi-isomorphism equivalence class of complexes
of bundles $(\mathcal{E}_i, d_{\mathcal{E}})$, whose charge is
\begin{equation}
\label{eq-QE }
Q_{{\mathcal{E}}^{*}}  = \sum_i (-)^i \ch(\mathcal{E_i}) \sqrt{\td( M)}.
\end{equation}
In the large volume limit of $ M$, its central charge is 
\begin{equation}
\label{centralchargeBbrane}
Z_{\mathcal{E}^{*}}= \int_ M{\exp^{-(B+iJ)( M)}}\ch(\mathcal{E}^{*}) \sqrt{\td( M)},
\end{equation}
which explicitly depends on the Kaehler moduli $(B+iJ)( M)$. 

In this work, we will restrict attention to the large volume limit,
loosely speaking $J >> 1$ (in string units) and $|F|,|B| << J$.  In
this limit, one expects the $\Pi$-stable objects to be $\mu$-stable
coherent sheaves which support solutions of the hermitian Yang-Mills
equations, possibly with mild singularities.  In \cite{BS}, it is
shown that such solutions exist for $\mu$-stable reflexive sheaves,
which motivates the precise form of our Conjecture 1.1.  Regarding
Conjecture 1.2, on surfaces such sheaves are necessarily locally free.

Let us discuss a few of the known stringy corrections to this limit.
Some of these can be taken into account by studying a generalization
of the hermitian Yang-Mills equations, the MMMS equations \cite{MMMS}.
It was shown by Leung \cite{Leung} that the solvability of these
equations follows from a deformed notion of stability.  This still
leaves further stringy corrections arising from world-sheet
instantons.  These can be best understood using the mirror symmetry to
the IIb string on $\tilde{M}$.  Now applying mirror symmetry directly
leads us to the problem of classifying special Lagrangian manifolds,
which at present appears more difficult than the original problem.
However, an indirect approach is to use mirror symmetry to motivate a
deformation of $\mu$-stability for bundles called $\Pi$-stability
\cite{Douglas:2000gi,DoA,BrS}, 
which is believed to incorporate all the stringy
corrections mentioned earlier.  A complete description along these
lines requires considering not just bundles, but arbitrary objects in
the derived category of coherent sheaves.  In the one case of a
Calabi-Yau threefold which has been completely analyzed at present
\cite{Br}, the total space of the line bundle
$\mathcal{O}_{\mathbb{P}^2}(-3)$, the set of $\Pi$-stable objects is
actually simpler in the stringy regime (say near the orbifold point)
than at large volume.  The set of stable objects varies continuously
with $[g]$, so one might be able to start with this simpler
description and evolve it up to large volume, providing a new approach
to this problem.  Doing this would also provide a precise definition
of the class of objects which can be used in this limit.

We should also discuss the precise singularities allowed in our
sheaves.  In general, string theory even allows other objects, such as
noncommutative analogs of bundles.  We do not know of a complete and
precise mathematical definition of the class of allowed objects, other
than the implicit definition provided by $\Pi$-stability.  It appears
to be different in the different contexts in which connections appear
in the large volume limit (heterotic string, D-branes, and at
singularities).  In the case at hand of B-type D-branes, it includes
the zero size limit of Yang-Mills instantons, and also includes
certain ``rank one'' or noncommutative instantons.  These correspond
to torsion-free coherent sheaves, for example the ideal sheaf of a
curve.

In this paper, we sidestep such questions by studying the predictions
of the attractor conjecture in the large volume limit, and accepting
the claim that these correspond to singular solutions of hermitian
Yang-Mills.  The conjecture itself is more general and it would be
interesting to compare it with the classification of these more
general objects.

To conclude the subsection, we will give an explicit map between
$H^{\mathrm{even}}( M,\mathbb{Q})$ and $H_3(\tilde{ M},\mathbb{Q})$. In
particular, we express the Chern character of any $\mu$-stable vector
bundle $V$ on $ M$ in terms of charges of special Lagrangian 3-cycles
on $\tilde{ M}$.

Consider a vector bundle $V$ which is stable with respect to some
Kaehler class $\omega$. Hence, in the large volume limit, $V$
corresponds to a BPS-state on $ M$. Using the coordinates
$\{\tilde{s}^{\alpha}\}$ on $\mathcal{M}_{Kaehler}( M)$ introduced in
\ref{eq-isomorphism(s) }, that is
\begin{equation}
B+iJ= \tilde{s}^{'\alpha} B_{\alpha} + i\tilde{s}^{''\alpha}J_{\alpha},
\end{equation}
where $\{B_{\alpha}\}$ denote a basis of $H^2( M,\mathbb{R})$ and $\{J_{\alpha}\}$ a framing of the Kaehler cone of $M$ (that is $\{J_{\alpha}\} \in H^2(B,\mathbb{Z})$ and $\{J_{\alpha}\}$ are in the closure of the Kaehler cone), we can compute the central charge $\mathcal{Z}_V(H)$ of the vector bundle $V$ at the point
\begin{equation}
\tilde{s}^{'\alpha}=0,\;\;\;iH\equiv \omega  =\tilde{s}^{''\alpha}J_{\alpha}
\end{equation}
Expanding  \ref{eq-QE }, we find for the charge $Q_V$ of the bundle $V$
\begin{equation}\label{eq-Z_V}
\begin{split}
           Q_V  = &  \ch(V) \sqrt{\td( M)} \\
               = & rk(V)+ c_1(V)+\left( ch_2(V)+\frac{rk(V)}{24}c_2( M)\right)+\left(ch_3(V)+\frac{1}{24} c_1(V)c_2( M)\right) \\
               = &\oplus_{i=0}^3 Q^{2i},\;\;\;\;Q^i \in H^i( M,{\mathbb{Q}}).\end{split}
\end{equation}
Using the expression \ref{centralchargeBbrane}, we find for the central charge of $V$ at the point $H$ 
\begin{equation}\label{eq-centralchargeofVatH}
\mathcal{Z}_V(H)=\frac{1}{6}H^3 Q^0 +\frac{1}{2}H^2 Q^2 + H Q^4 + Q^6.
\end{equation}
The central charge of the corresponding Lagrangian 3-cycle $\gamma_V$ on $\tilde{ M}$ is given by
\begin{equation}
\label{equalcentralcharge}
\mathcal{Z}_{\gamma_V}(H)\equiv \mathcal{Z}_V(H).\footnote{Up to a possible normalization.}
\end{equation}
More concretely, if we expand the charge $\hat{\gamma}_V$ in the symplectic basis of $H^3(\tilde{ M},\mathbb{Z})$ introduced in \ref{eq-symlectic}, that is
\begin{equation}
\hat{\gamma}_V=  ( p^I \hat{\alpha}_I - q_I \hat{\beta}^I),
\end{equation}
we find for $\mathcal{Z}_{\gamma_V}$ at the point ${s}^{\alpha}$ in the large complex structure limit o<f moduli space $\mathcal{M}_{Complex}{(\tilde{ M})}$
\begin{equation}\label{eq-centralchargeinpqtwo}
\begin{split}
Z_{\gamma_V}({s}) = & \int_{\gamma} \Omega \\
  = &    q_I \int_{\alpha^I}\Omega - p^I \int_{\beta_I} \Omega \\
 =&  X^0(q_0 +q_{\alpha}{s}^{\alpha} -p^0\mathcal{F}_0(s)-p^{\alpha}\mathcal{F}_{\alpha}(s)).
\end{split}
\end{equation}
It follows from local isomorphism \ref{eq-isomorphism(s) } that we need to evaluate the central charge $Z_{\gamma_V}({s})$ at
\begin{equation}
{s}^{'\alpha}=0,\;\;\;iH  = {s}^{''\alpha}J_{\alpha}
\end{equation}
in oder to compute $Z_{\gamma_V}(H)$.

To begin with, we compute the  in-homogenous pre-potential  $\mathcal{F}(H)$. We find for the coefficients $D_{\alpha\beta\gamma}$ and $B_{\alpha}$
\begin{equation}
\frac{1}{3!}D_{\alpha\beta\gamma}s^{\alpha}s^{\beta}s^{\gamma}=\frac{1}{3!}H^3
\end{equation}
and
\begin{equation}
 B_{\alpha} s^{\alpha}= -\frac{1}{24} i  s^{''\alpha}\int_ M J_{\alpha} \wedge c_2( M)  =\frac{1}{24}H c_2( M) 
\end{equation}
The integral matrix $A_{\alpha\beta}$ has no known topological interpretation and is usually fixed by monodromy transformations. We will consider A as a bilinear map $H^2(M)\times H^2(M) \to \mathbb{R}$ an leave it undetermined. Then we find
for the inhomogenous prepotential
\begin{equation}
\mathcal{F}(H) = \frac{1}{3!}H^3 +\frac{1}{2}AH^2+ \frac{1}{24}c_2( M)H
\end{equation}  
Using $\frac{\partial H}{\partial s^{\alpha}}= -J_{\alpha}$ we find 
\begin{equation}
\mathcal{F}_{\alpha} (H)= \partial_{\alpha} \mathcal{F} = -\frac{1}{2}  H^2  J_{\alpha} - AHJ_{\alpha}-\frac{1}{24} c_2( M) J_{\alpha}
\end{equation}
and
\begin{equation}
\mathcal{F}_0(H)=-\frac{1}{6}H^3+\frac{1}{24}c_2( M)H.
\end{equation}
We denote the dual basis of $\{J_{\alpha}\}$ by $\{C^{\alpha}\}$ and find for central charge $\mathcal{Z}_{\gamma_V}(H)$ 
\begin{equation}
\label{eq-Z_gamma}
\mathcal{Z}_{\gamma_V}(H)= X^0\left(q_0-q_{\alpha}C^{\alpha}H +p^0\left(\frac{1}{6}H^3-\frac{1}{24}c_2( M)H\right)+p^{\alpha}J_{\alpha}\left(\frac{1}{2}H^2+AH+\frac{1}{24}c_2( M)\right)\right)
\end{equation}
We expand both sides of \ref{equalcentralcharge} in powers of $H$ to derive 
the map  $H^{2*}( M,\mathbb{Q}) \to H^3(\tilde{ M},\mathbb{Q})$. In particular, we  express the vector $(q_I,p^J)$ in terms of the Chern classes of V. 
To begin with consider the term proportional to $H^3$. We find
\begin{equation}
 X^0 \frac{1}{6}p^0H^3 \sim \frac{1}{6}H^3 Q^0=\frac{1}{6}H^3 rk(V).
\end{equation}
We will assume
\begin{equation}
\label{eq-rankV}
p^0= rk(V),
\end{equation}
which fixes the normalization. After comparing the terms proportional to $H^i$ for $i=2,1,0$ we find
\begin{equation}\label{eq-chernV}
\begin{split}
  p^{\alpha}  J_{\alpha} &= c_1(V)\\
  q_{\alpha}    C^{\alpha}& =-\left(ch_2(V)+\frac{rk(V)}{12}c_2( M)\right)+c_1(V)A \\
       q_0 &= ch_3(V)      
\end{split} 
\end{equation}
Note that we have computed that map in the neighborhood of the
specific point $H \in \mathcal{M}_{Kaehler}( M)\approx
\mathcal{M}_{complex}(\tilde{ M})$. However, since the map is
topological, henceforth we will assume it holds everywhere in a
neighborhood of the large volume limit of $ M$ and $\tilde{ M}$.

\section{ Attractor equations}\label{s-solution}

The upshot of the preceding section is that we can phrase the IIa
attractor problem in the same terms as the IIb problem, by introducing
the following objects in $H^{2*}(M,\bbc)$,
\be \label{eq-defOmega}
\hat\Omega \equiv e^{B+iJ} ,
\ee
which spans the subspace of $H^{2*}$ mirror to $H^{(3,0)}(\tilde  M)$, and
\be
\gamma(V) = \Tr e^F \sqrt{{\rm Td~ M}}
\ee
which is mirror to the three-cycle wrapped by the D$3$-brane.
This neglects all stringy $\alpha'$ corrections; we will discuss these
later.

A IIa attractor point is then a point in complexified K\"ahler moduli
space for which $\gamma(V)$ is contained in the mirror to
$H^{(3,0)}\oplus H^{(0,3)}$, in other words
\be\label{eq-attractortwo}
\gamma(V) = \Re {\bar C} \hat\Omega
\ee
for some complex number $\bar C$. 
Explicitly, we expand 
$$
\hat\Omega = 1 + (B+iJ) + \frac{1}{2}(B+iJ)^2 + \frac{1}{6}(B+iJ)^3
$$
and 
$$
\gamma(V) = r + c_1 +ch_2+\frac{r}{24}c_2( M) + ch_3+\frac{1}{24}c_1c_2( M) .
$$
We will consider  \eq{attractor}  in the case for non-vanishing rank and vanishing rank separately.

\subsection{$V$ is a vector bundle on $M$}
We can solve the zero-form term in \eq{attractor}  by writing
$$
\bar C=r(1-i\xi)
$$
with a free parameter $\xi$.
The two-form term is then
$$
c_1 = \Re r(1-i\xi)(B+iJ) = r( B + \xi J)
$$
which is solved by writing
\be
\label{eq-solveB}
B = \frac{c_1}{r} - \xi~ J ,
\ee
thus determining the $n+1$ real parameters $B$ and $\Re C$ in terms of
the $n+1$ real parameters $r$ and $c_1$.

To find \ATT, the region in $H^{2*}(M)$ which supports attractor
points, we can now simply choose $r$ and $c_1$, and vary the possible
choices of $J$ and $\xi$, finding the corresponding $ch_2$ and $ch_3$.
Furthermore, in the large volume limit, the $c_1$ dependence is
essentially trivial, as using \eq{solveB} we can write
$$
\hat \Omega = e^B\ e^{iJ} = e^{c_1/r}\ e^{(i-\xi)J} .
$$
Since $e^{c_1/r}$ is real, we can solve \eq{attractor} with $c_1=0$, and then
multiply the result by $e^{c_1/r}$ to get the most general point in \ATT.

Thus we now set $c_1=0$, that is, we replace the chern classes of $V$ with the Chern classes of $V\otimes \mathcal{O}_ M(-c_1/r)$, and find
\be
\gamma(V) = \Re r(1-i\xi)\ e^{(i-\xi)}J  
\ee
which can be simplified to
\be \label{eq-attractor-sol}
r -\left( c_2 -\frac{r}{24}c_2( M)\right)+ \frac{1}{2}c_3 
= r - \frac{r}{2}(1+\xi^2)J^2 + \frac{r\xi(1+\xi^2)}{3}J^3 ,
\ee
where we used $c_2=-ch_2$ and $ch_3=\frac{1}{2}c_3$ in the case $c_1=0$.

These are the desired $n+1$ equations for the $n+1$ parameters $J$ and $\xi$.
Clearly the main effort in solving them is to find a real two-form $\tH$ satisfying
$$
(\tH)^2 = \frac{1}{r}\left(c_2-\frac{r}{24}c_2( M)\right) .
$$
We choose this normalization of $\tH$
to have invariance under the rescaling $ch\rightarrow N~ch$.
Then, we can write
\be\begin{split}
J &= \lambda \tH \\
\lambda^2 &= \frac{2}{1+\xi^2} \\ 
c_3 &= \frac{2r\xi(1+\xi^2)}{3}\lambda^3\tH^3  \\
&= \frac{2^{5/2}r}{3}\cdot \tH^3 \cdot \frac{\xi}{(1+\xi^2)^{1/2}} .
\end{split}\ee
We note that $\lambda$ and thus $J$ are also
invariant under the rescaling $ch\rightarrow N~ch$.

What are the consistency conditions on such a solution?  
We define $\ATT$ to be the largest 
set of solutions we might consider, by requiring
the class $\tH$ to lie in the closure of the K\"ahler cone.
This implies that
$$
\left(2rc_2(V)-\frac{r^2}{12}c_2(M)\right) \cdot J \geqslant 0 .
$$
for all $J$. 
Since for Calabi-Yau threefolds $c_2(M)\cdot J \geqslant 0$, any point
in $\ATT$ will satisfy the Bogomolov bound, for any $J$.

In deciding which points in $\ATT$ correspond to stable bundles, we
must discuss the stringy corrections to the problem.  First, there are
corrections to \eq{defOmega}.  These are in principle entirely determined
by mirror symmmetry and the considerations of the previous section;
we have
\be \label{eq-defOmegacor}
\hat\Omega \equiv e^{B+iJ} +i \frac{\zeta(3)\chi(M)}{(2\pi)^3} \omega^3
 + {\mathcal O}\left(e^{-J}\right) .
\ee
The form $\omega^3$ is a six-form whose integral on $M$ is $\alpha'^3$,
{\it i.e.} $1$ in ``string units.''  Similarly, the exponentially small
terms are functions of volumes of cycles measured in string units.

Since the original problem of solving the Hermitian Yang-Mills
equations or finding $\mu$-stable bundles did not require defining the
string unit of length, it would appear that these corrections are not
directly relevant, and we will neglect them in the following.
However, before proceeding, we should explain what role they play in
the full problem arising in string theory.  First, there are $\alpha'$
corrections to the HYM equations and the stability condition as well,
as discussed at length in \cite{Asp}.  Second, even if we start at
large volume (meaning $J>>1$ in string units, so that world-sheet
instanton corrections are negligible), the attractor point might be at
small volume.  We would then expect the BPS state to be stable at the
attractor point, where stringy corrections might be important.

For both reasons, a full treatment of the problem would involve
comparing the exact attractor points to the set of $\Pi$-stable
objects.  However, since (as we will see later) there are unresolved
issues even before we reach this point, we will leave this for future
work.

Given $c_2$ and $\tH$, the remaining step is to vary $\xi$ and see
what values of $c_3$ can be attained.  Since the quantity
$\xi/{(r^2+\xi^2)^{1/2}}$ takes all values between $-1$ and $1$, we
conclude that all $c_3$ can be attained which satisfy
\be\label{eq-cthreeboundtwo}
 |c_3| \le \frac{2^{5/2}}{3} r\cdot \tH^3.
\ee
Note that this bound implies that 
\be\label{eq-cthreeboundthree}
 |c_3| \le \frac{2^{5/2}}{3} r \left(\left(\frac{c_2}{r}-
 \frac{1}{12}c_2(X)\right)\cdot w\right)^{3/2}(w^3)^{-1/2}.
\ee
for any ample class $w$. This follows from \eq{ample}.

Finally, if $\xi$ becomes too large, $J$ will no longer be large.  It
is not completely clear whether this is a problem or what lower bound
to take, but if we require $J>1$, then $\lambda^2 c_2 > r$ so we need
$\xi^2 < c_2/r$.  This will bring the bound in \eq{cthreebound} down
by a rough factor $1-r/2c_2$, unimportant for $c_2>>r$ as typical in
the heterotic string application, but possibly important in general.
Thus the funny coefficient $2^{3/2}/3$ does not have as much
fundamental significance as one might think.

\subsection{$V$ is a sheaf on $M$ supported on a smooth hypersurface}

We consider the case that $V$ is supported on a smooth surface $i: D
\subset M$. More precisely, we consider $V$ to be stable bundle on
D. This corresponds to case that the rank of $V$ as a sheaf on $ M$
vanishes. To begin with, we compute the charges of $V$ in
$H^{2*}(M,\bbc)$, given its Chern classes on $D$. This can be
accomplished using the formula of Grothendieck-Riemann-Roch. For
notational simplicity we will denote the Chern character of $V$ on $D$
by $ch$ and on $M$ by $ch \ i_{*}V$. GRR states that
\begin{equation}
\ch(i_{*}\mathcal{F})=i_{*}(\ch(\mathcal{F})\td(D))\td( M)^{-1},
\end{equation}
Using
\begin{equation}
c_1(D)=-D^2,\;\;\;c_2(D)=D^3+c_2( M)D
\end{equation}
we find
\begin{equation}
\begin{split}
r(i_{*}\mathcal{V})&=0\\
ch_1(i_{*}\mathcal{V})&= rD\\
ch_2(i_{*}\mathcal{V})&=-\frac{rD^2}{2}+c_1\\
ch_3(i_{*}\mathcal{V})&=\frac{r)D^3}{6}+ch_2-\frac{(D^2) \cdot  c_1}{2}.\\
\end{split}
\end{equation}
and 
\begin{equation}
\gamma(i_{*}V)=rD-\left(\frac{rD^2}{2}-c_1\right)+\left(\frac{rD^3}{8}+ch_2-\frac{(D^2) \cdot  c_1}{2}+\frac{r}{24}c_2(D)\right)
\end{equation}
The zero and two form term in \eq{attractor} are solved by
\begin{equation}
\bar{C}=-i\xi,\;\;J=\frac{r}{\xi}D
\end{equation}
where $\xi$ is a free real parameter. To solve the four form term we have to find a $B$ such that
\begin{equation}
\label{eq-W }
c_1(V)=rD\left(B+\frac{D}{2}\right).
\end{equation}
This equation simply implies that $c_1$ is given as the restriction of a class in $H^2$. 
We will denote this lift  by $\tilde{c_1}$. The strong Lefschetz theorem implies that the intersection
\begin{equation}
D: H^{(1,1)} \to H^{(2,2)}
\end{equation}
is an isomorphism. We can invert \eq{W } and find
\begin{equation}
B=\frac{1}{r}{\tilde{c_1}}-\frac{D}{2}.
\end{equation}
Finially, the six form part of \eq{solveB} implies
\begin{equation}
\label{bogsurfaces}
2rc_2-(r-1)c_1^2-\frac{r^2}{12}c_2(D)=\frac{r^2 D^3}{3 \xi^2}.
\end{equation}

\section{Multi-center attractor solutions}

Our main result can be summarized in the conjecture
$$
\ATT \subset \PCHu,
$$ 
where  $\ATT \subset H^{2*}(M,\bbr)$ is the set of attractor points,
$$
\ATT = \{ \Re {\bar C} e^{B+iJ} | C\in\bbc, B\in H^2(M,\bbr), J\in \KC(M)\}
$$
(here $\KC\subset H^2(M,\bbr)$ is the K\"ahler cone), and
$\PCHu$ is the projectivized set of Mukai vectors \eq{mukai} for which
a $\mu$-stable bundle exists for some $J$.

This conjecture passes the only definitive results we know on
$\CHu$, namely $\ATT$ satisfies the Bogomolov bound, and realizes a
bound on $c_3$ as required by Maruyama-Langer.  
In fact, we obtained a fairly explicit sufficient condition
on $c_3$, \eq{cthreebound}.  


How similar is $\ATT$ to $\bar{CH}$ ?  For K3 surfaces, they are the
same.  However, for threefolds it is easy to find stable bundles whose
Chern character is not in $\ATT$. The tangent bundle on the quintic
hypersurface $TQ$ is already an example: we have $c_3(TQ)=-200$, while
from $\tilde{H}^2=\frac{70}{24}H^2_Q$ follows that
$\tilde{H}^3\simeq5$ which violates the bound.
Another set of examples are given by Maruyama's construction. 
Recall from subsection \ref{ss-suff} that these can realize
\begin{equation}
c_2\sim n,\;\;\;c_3\sim n^2,
\end{equation} 
and therefore at fixed $r$ and $H$, for large $n$ our bound will be
violated. Thus, we can only claim that it is a sufficient condition,
which is all that our arguments guaranteed so far.

\subsection{Bound states}

What are the physical states corresponding to stable bundles violating
our bounds? So far, we have only considered BPS states corresponding
to spherical symmetric black hole solutions in four dimensions.  A
more general class of solutions is the bound states of spherical
symmetric solutions, as discussed by Denef \cite{De}.  Such a solution
will have a charge which is the sum of charges of two or more
constituents, each of which should lie in $\ATT$.  
Because of these solutions, the condition
$Z(\gamma,{B+iJ})=0$ does not rule out the existence of a BPS state of
charge $\gamma$, as we now explain.

Not all sums of elements of $\ATT$ correspond to multi-center solutions.
As shown in 
\cite{De,Denef:2001xn,Denef:2001ix,Denef:2002ru,Bates:2003vx},
a multi-center solution corresponds to a ``split attractor'' flow.
Mathematically, this is an embedding of an oriented tree into moduli space,
such that the edges are gradient flows, the endpoints are attractor points,
and the vertices are additive decompositions of the charge vector,
$$
\gamma = \gamma' + \gamma'',
$$
satisfying the following
two consistency conditions.  First, the intersection product \eq{intersect}
must be non-zero, $\chi(\gamma',\gamma'')\ne 0$.  Second, the vertex must
embed into a point at which $Z(\gamma')\ne 0$, $Z(\gamma'')\ne 0$ and
$Im (Z(\gamma')/Z(\gamma'')) = 0$; in other words
on a line of variation of stability.

The start point (call it $t_0$) of the split attractor flow (as for
other attractor flows) is physically the value of the moduli ``in the
asymptotic region,'' so mathematically the existence or non-existence
of a flow implies stability at $t_0$, either $\mu$-stability for $t_0
>> \gamma$ (the large volume limit), or $\Pi$-stability in general.

In general, the analysis of the existence of these flows is rather
intricate, and in the references has been done by explicit numerical
integration of the attractor equations and search over candidate
additive decompositions.  However there is one simple consequence if
an attractor point is a regular zero.
It is clear from
the definition of split attractor flow that one can only exist if
the start point $t_0$ is different from the attractor point $B+iJ$.
Therefore, if $Z(\gamma,{B+iJ})=0$ at a point $B+iJ$
in the interior of moduli space, we can infer that a BPS state
of charge $\gamma$ {\it cannot} be stable at that $B+iJ$.

In general, this is a statement about $\Pi$-stability.
It would be of direct interest for our original goal of characterizing 
$\mu$-stable bundles, if the zero was attained in the large volume
limit, $J >> |\ch(V)|$ (taking say the $l_1$ norm).
However, it is easy to see that this will not happen, using the
fact that large volume formula \eq{centralchargeBbrane} is polynomial in
$\ch(V)$, and elementary bounds on the locations of zeroes of polynomials.
Still, it is interesting for example that $\mu$-stable bundles which violate
\eq{cthreeboundtwo} will have lines of marginal stability at some
$J\sim\ch V$, analogous to those for high degree hypersurfaces
\cite{DoA,De}.

To further illustrate the nature of the split attractor condition, let
us now consider a simpler set of necessary conditions for the
existence of a split attractor flow.

One simple criterion for when two states of charges $\gamma^{'}$ and
$\gamma^{''}$ in $H^{3}(\tilde{M})$ and central charge $Z^{'}$ and
$Z^{''}$ can form a bound state is that
\begin{equation}
\label{eq-bps-bound}
<\gamma^{'},\gamma^{''}>Im(Z^{'}\bar{Z}^{''})\geqslant 0.
\end{equation}
This condition follows from the analysis of the gradient flow for
$Z{'}, Z{''}$ and $Z$. To see this, if we assume that the gradient
flow for a BPS state of charge $\gamma$ (parameterized by the time
$\tau \in (0,\infty)$) crosses a line of variation of stability and
decomposes into two states of charge $\gamma'$ and $\gamma''$
respectively, then the time $\tau_{vs} $ can be explicitly computed by
\begin{equation}
\label{ }
\tau_{vs}=
 \frac{2Im(Z'\bar{Z}'')}{|Z|<\gamma',\gamma''>}\bigg|_{\tau=0},
\end{equation}
hence \ref{eq-bps-bound} follows.

One expects that on the bundle side, such a bound state will correspond to an
extension of bundles,
\begin{equation}
\label{eq-ext}
0 \to V^{'} \to V \to V^{''} \to 0,
\end{equation}
which can exist if $\Ext^1(V^{''},V^{'})\ne 0$ and $\mu(V^{'})<\mu(V^{''})$.
For BPS states, there is not such a clear distinction between positive and
negative rank, and one might also find bound states which correspond
to sub or quotient bundles,
\begin{equation}
\label{eq-sub}
V^{'}\to V^{"}\to W,\;\;or\;\;W \to V^{'}\to V^{"} .
\end{equation}
These will exist if there is a non-trivial map
$\alpha\in \Hom(V^{'}, V^{"})$
with vanishing kernel or cokernel.

On the bundle side, these are not easy conditions to check in general,
and one cannot expect to find a simple sufficient condition like
\eq{bps-bound}.  However, the hope would again be that in an
asymptotic limit of large Chern character, a simple condition emerges.
Thus, let us use the correspondence between special Lagrangian
three-cycles and stable bundles \eq{rankV}, to translate
\eq{bps-bound} into a statement about bundles.

We find for two bundles $V^{'}$ and
$V^{''}$ with Chern character $\ch^{'}$ and $ \ch^{''}$
\begin{equation}
<\gamma^{'},\gamma^{''}>=\chi(hom(V^{''},V^{'})).
\end{equation}
Expanding $Im(Z^{'}\bar{Z}^{''})$ in orders of $J$ (we will neglect $B$) gives
\begin{eqnarray}
\frac{1}{r^{'}r{''}}Im(Z^{'}\bar{Z}^{''})&=& -\frac{J^3}{12}\left(\mu^{''}-\mu^{'}\right)
  +  \frac{J^3}{6}\left(\frac{\ch_3^{''}}{r^{''}}-\frac{\ch_3^{'}}{r^{'}}\right) \\
 & + & \frac{J}{2}\left(\left (\frac{\ch_2^{'}}{r^{'}}+\frac{1}{24}c_2(X)\right)\mu^{''}-                           
 \left( \frac{\ch_2^{''}}{r^{''}}+\frac{1}{24}c_2(X)\right)\mu^{'} \right) \\
 &-& {J}\left( \left(\frac{\ch_2^{'}}{r^{'}}+\frac{1}{24}c_2(X)\right)-\left( \frac{\ch_2^{''}}{r^{''}}+\frac{1}{24}c_2(X)\right) \right)
\end{eqnarray}
where $\mu^{'}$ and $\mu^{''}$ denotes the slope of $V^{'}$ and $V^{''}$ respectively. Therefore, in the large volume limit, the above conditions reads as
\begin{equation}
\chi(hom(V^{'},V^{''}))(\mu^{''}-\mu^{'})\geqslant 0.
\end{equation}
Let us assume $\mu^{'}< \mu^{''}$. Then our condition translates into
the positivity of the Euler characteristic $\chi(hom(V^{'},V^{''}))$.
This requires either $Hom(V^{'},V^{''})> 0$ or $
Ext^2(V^{'},V^{''})=Ext^1(V^{''},V^{'})> 0.  $

Thus, the condition \eq{bps-bound} almost corresponds to the condition
required that one of a stable extension, sub or quotient bundle can be
constructed from the pair $V^{'}$ and $V^{''}$, missing only the condition
of vanishing kernel or cokernel of $\alpha$.  Optimistically assuming that
for large Chern character this condition becomes generic, it is reasonable
to expect the closure of $\ATT$ under this construction to become a better
approximation to $\CH(J)$.

Thus, we define the ``$J$-closure'' $S_J$ of a subset $S\subset\cH$ to
be the smallest set containing $S$ and all sums $\gamma'+\gamma''$
made from pairs satisfying \eq{bps-bound}.

Let us compare $\ATT_J$ to $\CH(J)$.  We begin with a pair of line
bundles $\om(p J)$ and $\om(q J)$ for some ample $J$; for $p>q$ one
expects to construct subbundles of the form
$$
0 \rightarrow E \rightarrow \om(q J)^m \rightarrow \om(p J)^n \rightarrow 0
$$
for various $m$ and $n$.
Since $\ATT$ and the condition \eq{bps-bound} are
homogeneous, we will not get a strong constraint on $m$ and $n$; indeed
this pair satisfies \eq{bps-bound} for $n>m$.  Thus, let us take
$m=p$ and $n=q$, to get $c_1(E)=0$.

We find that $\ATT_J$ contains the Chern characters
$$
\ch(E) = (p-q) + \frac{pq(q-p)}{2} J^2 + \frac{pq(q^2-p^2)}{6} J^3 .
$$
The ratio $c_3^2 r/c_2^3 \sim (p+q)^2/pq$ grows without bound,
so these certainly include points not in $\ATT$.  

On the other hand, they appear consistent with (an assumed homogeneous
version of) Maruyama-Langer, as we have bounded
$$
|c_3| r/c_2^2 \sim (p+q)/pq .
$$

Unfortunately, if we proceed to consider bound states of a pair of these
bundles, we find that these can attain arbitrarily large $c_3$ at fixed
$c_2$, so they satisfy no interesting bounds at all.  
Now as we explained, we have only considered a subset of the necessary
conditions for a split attractor flow, so there is no contradiction at
this point; rather we conclude that a more detailed analysis is required.

Such an analysis appears rather non-trivial and thus let us outline
some of the possible outcomes as a guide for future work.

One 
immediate mathematical question is whether the actual bounds on
Chern characters of $\mu$-stable bundles are homogeneous or not.
If not, we cannot expect to duplicate the bound by this analysis.

It is not at present clear to us whether the bounds \eq{gbound}
discussed in \cite{La2} are homogeneous or not.  If they are, then a
necessary condition on the third Chern class (for bundles with
vanishing first Chern class) could take the form
\begin{equation}
\label{eq-guess}
\begin{split}
|c_3|&\leqslant \frac{2\zeta(3)|\chi(M)|}{(2\pi)^3}r \left(\frac{c_2 \cdot J}{r}\right)^{1/2}(J^3)^{-1/6}\\
&+\frac{2^{5/2}}{3}r\left(\left(\frac{c_2}{r}-\frac{1}{12}c_2(X)\right) \cdot J\right)^{3/2}(J^3)^{-1/2}\\
&+ c r\left(\frac{c_2 \cdot  J}{r}\right)^{2}(J^3)^{-2/3}\\.
\end{split}
\end{equation}
for any ample class $J $ with respect the bundle is stable. The first
term is the correction which comes from the term of order $\alpha^3$
in \ref{eq-defOmegacor}, the second term is our usual bound for
elements in $\ATT$ and the third term (with undetermined constant $c$)
we have included to obtain consistency with known constructions, in
particular, with our examples from the monad construction.

To summarize, our attempts to enlarge $\ATT$ to include bound states
as a candidate description of all of $\CHu$ remain inconclusive.  One
possible explanation is that the actual mathematical bound on $c_3$ is
not homogeneous, in such a way that the set $\PCH$ simply contains all
charges consistent with Bogomolov.  Our simplified criterion for bound
states reproduces this, leading to a consistent if uninteresting
picture.  An alternative possibility is that a homogeneous necessary
condition exists, such as \eq{guess}, and that a refined version of
the bound state condition would make contact with this.  Finally, it
is conceivable that, even if the correct mathematical bound is not
homogeneous, incorporating higher genus corrections to the attractor
conditions (as discussed in
\cite{LopesCardoso:1998wt,Ooguri:2004zv})
which do not respect
homogeneity might provide a better description of $\CH$.

Finally, we should note that the definition of $\ATT$ given in section
4.1 makes sense in arbitrary dimensions, not just $d\le 3$, and it would
be interesting to test the analogous conjecture there.  While the
physical arguments do not directly apply to this case, many elements of
the physical discussion (notably the existence of boundary states and
the central charge formulas) do generalize, at least to Calabi-Yau
manifolds of arbitrary dimension.  One cannot expect the conjecture to
give a very good picture of $\CH$ for $d\ge 4$; for example we have
$\dim_{\bbr}\ATT = 2+2b^{1,1}<\dim\CH$, but the claim that
$\ATT\subset\CH$ still appears reasonable.

\vskip 0.2in
\noindent {\bf Acknowledgements}

We thank F. Bogomolov, F. Denef, R. Donagi, A. Langer, G. Moore and T. Pantev
for valuable discussions, and particularly F. Denef for a critical
reading of the manuscript and numerous suggestions.

This work was supported in part by DOE grant DE-FG02-96ER40959.

\vskip 0.2in
\noindent {\bf Appendix A}

On solid mathematics grounds, the only known general result is the
Bogomolov bound \eq{Bogomolov}. This bound however can be improved if
we restrict to certain algebraic surfaces.  First note that stable
bundles must be simple, that is
\begin{equation}
V\otimes V^{*}=End(V)\oplus \mathcal{O}_ D
\end{equation}
with
\begin{equation}
H^0(D,End(V))=0.
\end{equation}
where $D$ denotes an algebraic surface.
If we consider vector bundles with vanishing first Chern class stability implies $H^0(D,V)=0$.
Let $D$ be an $K3$ surface. Using Serre duality one finds $h^2( D,V)=h^0( D,V^{*})=0$. Computing the index for $V$ on $D$ on finds
\begin{equation}
\label{eq-Bogsur}
rc_2-\frac{r^2}{12}c_2(D)=-h^1(D,V)\geqslant 0.
\end{equation}
For bundles with non-vanishing first Chern class on can computes the index of  $V\otimes V^{*}$.  Since $V\otimes V^{*}$ is self-dual, $c_1(V\otimes V^{*})$ vanishes and one is in the case the case above. We find
\begin{equation}
2rc_2-(r-1)c_1^2-\frac{r^2}{12}c_2( D)=-\sum_i (-1)^{i}h^i(D,V\otimes V^{*}) \geqslant -2.
\end{equation}
Note that is this condition agrees with \eq{Yos}. Also note that the
bound will be saturated only by bundles exceptional bundles.  Lets us
know assume $D$ to be fano, that is, the anticanonical bundle is
ample. Well know examples of fano surfaces are $\mathbb{P}_2$, the
Hirzebruch surfaces $\mathbb{F}_r$ and the del Pezzo surfaces
$d\mathbb{P}_m,\;m=1,...,8$.  Consider stable bundles $V$ with
$c_1(V)=0$. Using Serre duality we find again $h^2( D,V)=h^0(
D,V^{*}\otimes K_ D)=0$, hence
\begin{equation}
rc_2(V)-\frac{r^2}{12}(c_2( D)+c_1( D)^2)\geqslant 0.
\end{equation}
In particular, for fano surfaces $D$ which admit an Kaehler Einstein metrics we have
\begin{equation}
c_2( D)+c_1( D)^2>0,
\end{equation}
hence the inequalities above  strengthen the Bogomolov bound.
So far we have not used the specific form of the ample class with respect $V$ is stable.  Consider a stable vector bundle $V$ with respect to an ample class $H$ such that $K_X  H<0$. (In the fano case this will hold for all ample classes). Then, for any surface, Maruyama  has shown that $Ext^2(V,V)=0$ and we find
\begin{equation}
\label{eq-indexsurfacegeneraltype }
2rc_2-(r-1)c_1^2-\frac{r^2}{12}\left(c_2(D)+c_1^2(D)\right)\geqslant -1.
\end{equation}
which will again be saturated for exceptional bundles.

Following our discussion on  the existence of stable vector bundles 
motivates the question whether the Bogomolov
bound can be improved in general. Note however, since we do not quite understand the split flow problem, the room for improvement may not be large. On the other hand, it is interesting for algebraic geometers to see, how much one can strengthen such an inequality.
Motivated by equation~\ref{bogsurfaces}, on might ask that on any simply connected surface $D$ with ample or trivial canonical bundle, the Chern classes of any stable  vector bundle with non-trivial moduli space and rank $r\geqslant2$ obey
\begin{equation}
\label{eq-conjBogomtwo}
2rc_2-(r-1)c_1^2-\frac{r^2}{12}c_2(D)\geqslant0.
\end{equation}
In  Appendix B we give an example which shows that this improvement  cannot be generalized to Calabi-Yau threefolds. 

\vskip 0.2in
\noindent {\bf Appendix B}
In this section we give an example of stable vector bundle $K$ of rank three on a generic quintic threefold $Q$ suggested by M. Jardim which shows that  the bound \ref{eq-conjBogomtwo} cannot be extended to Calabi-Yau threefolds.
Define $K$ as the kernel of the map $\beta$
$$
\xymatrix{
K \ar[r] &{\cal O}^{\oplus 4}_Q  \ar[r]^{\beta} &{\cal O}_Q(1)
}
$$
where $\beta$ is given by four generic global sections in ${\cal O}_Q(1)$ and ${\cal O}_Q(1)$ corresponds to the restriction of the hyperplane bundle ${\cal O}(H)_{\mathbb{P}^4}$. The Chern classes of $K$ are
$$
c_1(K)=-H,\;\;c_2(K)=H^2.
$$
Hence we find
$$
(2rc_2-(r-1)c_1^2)\cdot H=4H^3,\;\;\frac{r^2}{12}c_2(Q)\cdot H = \frac{15}{2}H^3
$$
which violates bound \ref{eq-conjBogomtwo}.
To proof stability we will show 
\begin{equation}
H^0(Q,K)=H^0(Q,\wedge^2 K) =0.
\end{equation}
Granted this fact and  that $Pic(Q)=\mathbb{Z}$, it is easy to see that $K$ is stable. Consider for example a possible destabilizing line bundle of $K$. It is of the form ${\cal O}_Q(p)$ form some $p \geqslant 0$. But this implies the existence of global sections of $K$, a contradiction.
From the defining sequence of $K$ it is clear that $H^0(X,K)$ vanishes for a generic map $\beta$. To see that $H^0(X,\wedge^2 K)=0$, consider
$$
h^0(X, \wedge^2 K)=h^0(X, K^* \otimes \det K)=h^3(X, K \otimes \det K^*)
$$
The vanishing of $h^3(X, K \otimes \det K^*)$ follows from the fact that ${\cal O}_Q(m)$ is ample for $m>0$.



\end{document}